\documentclass[12pt]{article}

\usepackage{algorithm, algpseudocode}
\usepackage{float}
\usepackage{tikz}
\usetikzlibrary{calc}
\usepackage{capt-of}
\usepackage{color}
\usepackage{fullpage}
\usepackage{amssymb}
\usepackage{amsmath}
\usepackage{amsthm}
\usepackage{amsfonts}
\usepackage{amscd}

\usepackage[colorlinks=true,
linkcolor=brown,
filecolor=brown,
citecolor=brown]{hyperref}

\newcommand{\seqnum}[1]{\href{https://oeis.org/#1}{\rm \underline{#1}}}

\begin{document}

\theoremstyle{plain}
\newtheorem{theorem}{Theorem}
\newtheorem{corollary}[theorem]{Corollary}
\newtheorem{lem}[theorem]{Lemma}
\newtheorem{proposition}[theorem]{Proposition}
\newtheorem{definition}[theorem]{Definition}

\begin{center}
\vskip 1cm{\LARGE\bf 
Counting self--dual monotone \\ Boolean functions
}
\vskip 1cm
\large
Bart\l{}omiej Pawelski\\ Institute of Informatics\\ University of Gda\'nsk\\ \href{mailto:bartlomiej.pawelski@ug.edu.pl}{\tt bartlomiej.pawelski@ug.edu.pl}

\vskip 1cm
\large
Andrzej Szepietowski\\
Institute of Informatics\\
University of Gda\'nsk\\
\href{mailto:andrzej.szepietowski@ug.edu.pl}{\tt andrzej.szepietowski@ug.edu.pl} \\
\end{center}

\vskip .2in

\begin{abstract}
Let $D_n$ denote the set of monotone Boolean functions with $n$ variables. Elements of $D_n$ can be represented as strings of bits of length $2^n$. Two elements of $D_0$ are represented as 0 and 1 and any element  $g\in D_n$, with $n>0$,  is represented as
a concatenation $g_0\cdot g_1$, where $g_0, g_1\in D_{n-1}$ and $g_0\le g_1$. 
For each $x\in D_n$, we have {\it dual} $x^*\in D_n $ which is obtained by reversing and negating
all bits.
An element $x\in D_n$ is {\it self-dual} if $x=x^*$.
 Let $\lambda_n$ denote the cardinality of the set of all self-dual monotone Boolean functions of $n$ variables. The value $\lambda_n$ is also known as the $n$-th Hosten-Morris number.
In this paper, we derive several algorithms for counting self-dual monotone Boolean functions and confirm the known result that $\lambda_9$ equals 423,295,099,074,735,261,880.
\end{abstract}
    
\section{Introduction}\label{sec:intro}

Let $B$ denote the set $\{0,1\}$ and $B^n$ the
set of $n$-element sequences of $B$.
A Boolean function with $n$ variables is any function
from $B^n$ into $B$. There is the order relation in $B$ (namely: $0\le 0$, $\;0\le1$, $\;1\le1$)
and the partial order in $B^n$:
for any two elements:
$x=(x_1,\dots,x_n)$, $\;y=(y_1,\dots,y_n)$ in $B^n$,
$x\le y$ if and only if $x_i\le y_i$ for all $1\le i\le n$.
A function $h:B^n\to B$ is monotone if $x\le y$ implies $h(x)\le h(y)$.
Let $D_n$ denote the set of monotone functions with $n$ variables
and let $d_n$ denote $|D_n|$.
We have the partial order in $D_n$ defined by:
$$g\le h\quad\hbox{if and only if}\quad g(x)\le h(x)\quad
\hbox{for all } x\in B^n.$$
We shall represent the elements of $D_n$ as strings of bits of
length $2^n$. Two elements of $D_0$ will be represented as 0 and 1.
Any element $g\in D_1$ can be represented as
the concatenation $g(0)\cdot g(1)$, where $g(0), g(1)\in D_0$ and $g(0)\le g(1)$. Hence, $D_1=\{00, 01, 11\}$. 
Each element of $g\in D_2$ is the concatenation (string) of four bits: $g(00)\cdot g(10)\cdot g(01)\cdot g(11)$ which can be represented as
a concatenation $g_0\cdot g_1$, where $g_0, g_1\in D_1$ and $g_0\le g_1$.
Hence, $D_2=\{0000, 0001, 0011, 0101, 0111, 1111\}$.
Similarly, any element of $g\in D_n$ can be represented as
a concatenation $g_0\cdot g_1$, where $g_0, g_1\in D_{n-1}$ and $g_0\le g_1$. 

For each $x\in D_n$, we have {\it dual} $x^*\in D_n $, which is obtained by reversing and negating all bits. For example, $1111^*=0000$ and $0001^*=0111$. An element $x\in D_n$ is {\it self-dual} if $x=x^*$. For example, $0101$ and $0011$ are self-duals in $D_2$.

 Let $\Lambda_n$ be the set of all self-dual monotone Boolean functions of $n$ variables, and let $\lambda_n$ denote the cardinality of this set. The value $\lambda_n$ is also known as the $n$-th Hosten-Morris number (\seqnum{A001206} in \textit{On-Line Encyclopedia of Integer Sequences}). The first attempt to solve the problem of determining the values of $\lambda_n$, as found in the literature, was made by Riviere \cite{riviere} in 1968, who determined all values up to $\lambda_5$. In 1972, Brouwer and Verbeek provided the values up to $\lambda_7$ \cite{brouwer7}. The value of $\lambda_8$ was determined by Mills and Mills \cite{mills} in 1978.

The most recent known term, $\lambda_9$, was obtained by Brouwer, Mills, Mills, and Verbeek \cite{brouwer} in 2013. The value of $\lambda_n$ also corresponds to the number of maximal intersecting families on an $n$-set \cite[Section 1]{brouwer}, as well as the number of non-dominated coteries on $n$ members \cite[Section 1]{bioch}.

In this paper, we derive several algorithms for counting self-dual monotone Boolean functions. We also confirm the result of \cite{brouwer} that $\lambda_9$ equals 423,295,099,074,735,261,880.

\begin{table}[ht]
\centering
\begin{tabular}[h]{|p{1.2cm}|l|}
\hline
 $n$ & $\lambda_n$\\
\hline
0 &  0 \\
1 &  1 \\
2 &  2 \\
3 &  4 \\
4 &  12 \\
5 &  81 \\
6 &  2,646 \\
7 &  1,422,564 \\
8 &  229,809,982,112 \\
9 &  423,295,099,074,735,261,880 \\

\hline
\end{tabular}
\caption{Known values of $\lambda_n$ (\seqnum{A001206})}
\label{tab:SD}
\end{table}

\section{Preeliminaries}

By $\top$ we denote the maximal element in $D_n$, that is $\top=(1\ldots 1)$, and by $\bot$ the minimal element in $D_n$, that is $\bot=(0\ldots 0)$. For two elements $x,y\in D_n$, by $x|y$ we denote the \textit{bitwise or}; and by $x\&y$ the \textit{bitwise and}. Furthermore,
let $re(x,y)$ denote $|\{z\in D_n : x\le z\le y\}|$. Note that $re(x,\top)=|\{z\in D_n : x\le z\}|$ and $re(\bot,y)=|\{z\in D_n : z\le y\}|$. 
For $x\in D_n$, by $\ell(x)$ we denote the number of ones in $x$, also known as its Hamming weight. For example,
$\ell(0000)=0$ and $\ell(0101)=2$.

\begin{lem}\label{L0}
For each $x,y\in D_n$, we have:
\begin{enumerate}
\item $x^{**}=x$
\item if $x\le y$ then $y^*\le x^*$
\item $(x|y)^*=x^*\&y^*$
\item $(x\&y)^*=x^*|y^*$
\end{enumerate}
\end{lem}

\subsection{Posets}
A {\it poset (partially ordered set)} $(S,\le)$ consists of a set $S$ (called the carrier) together with a binary relation (partial order) $\le$ which is reflexive, transitive, and antisymmetric.
For example, $B$, $B^n$, and $D_n$ are posets. 
Given two posets $(S,\le)$ and $(T,\le)$, a function $f:S \to T$ is {\it monotone}, if $x \le y$ implies $f(x)\le f(y)$.
By $T^S$ we denote the poset of all monotone functions from $S$ to $T$ with the partial order  defined by:
$$f\le g\quad\hbox{ if and only if}\quad f(x)\le g(x)\hbox{ for all }x\in S.$$
Notice that $D_n=B^{B^n}$ and $B^n=B^{A_n}$, where $A_n$ is the antichain with
the carrier $\{1,\ldots,n\}$.
In this paper we use the well known lemma:
\begin{lem}\label{L1} The poset $D_{n+k}$ is isomorphic to the poset $D_n^{B^k}$---the poset of monotone  functions from $B^k$ to $D_n$.
\end{lem}
\subsection{Permutations and equivalence relation}
Let $S_n$ denote the set of permutations on $\{1,\ldots,n\}$.
Every permutation $\pi\in S_n$ defines the permutation on $B^n$ by
$\pi(x)=x\circ\pi^{-1}$ (we treat each element $x\in B^n$ as a function
$x:\{1,\ldots,n\}\to \{0,1\}$).
The permutation $\pi$ also generates the permutation
on $D_n$. Namely, by $\pi(g)=g\circ \pi$.
By $\sim$ we denote an equivalence relation on $D_n$. Namely,
two functions $f, g\in D_n$ are {\it equivalent}, $f\sim g$, if there is a
permutation $\pi\in S_n$ such that $f=\pi(g)$.
For a function $f\in D_n$ its {\it equivalence class}
is the set $[f]=\{g\in D_n : g \sim f\}$. By $\gamma(f)$ we denote $|[f]|$.
For the class $[f]$, its {\it  representative} is its minimal element
(according to the total order induced on $D_n$ by the total order in integers). Sometimes, we identify the class $[f]\in R_n$ with
its representative and treat $[f]$ as an element in $D_n$.
By $R_n$ we denote the set of equivalence classes and by $r_n$ we denote the number of the equivalence classes; that is $r_n=|R_n|$.
\begin{lem}\label{L2}

\begin{enumerate}

\item For every element $x\in D_n$ and every permutation $\pi\in S_n$, we have
$\pi(x^*)=(\pi(x))^*$.        
\item If $x\in D_n$ is self-dual, then every equivalent $y\in[x]$ is self-dual.
\end{enumerate}
\end{lem}

\section{Counting functions from $B^2$ to $D_n$}
By Lemma~\ref{L1}, the poset $D_{n+2}$ is isomorphic to the poset $D_n^{B^2}$---the poset of monotone functions from $B^2=\{00,01,10,11\}$ to $D_n$. 
Consider a monotone function $H:B^2\to D_n$. It can be represented as the concatenation:
$$H(00)\cdot H(01)\cdot H(10) \cdot H(11)$$
and its dual as
$$H(11)^*\cdot H(10)^*\cdot H(01)^* \cdot H(00 )^*.$$
If $H$ is self-dual then $H(00)=H(11)^*$, $H(01)=H(10)^*$, $H(10)=H(01)^*$, and $H(11)=H(00)^*$.

\def\drawconnectionsa{\draw (0,0) -- (-2,2);
    \draw (0,0) -- (2,2);
    \draw (-2,2) -- (0,4);
    \draw (2,2) -- (0,4);}
    
\def\drawconnectionsb{\draw (0,0) -- (-2,2);
    \draw (0,0) -- (0,2);
    \draw (0,0) -- (2,2);
    \draw (-2,2) -- (-2,4);
    \draw (-2,2) -- (0,4);
    \draw (0,2) -- (-2,4);
    \draw (0,2) -- (2,4);
    \draw (2,2) -- (0,4);
    \draw (2,2) -- (2,4);
    \draw (-2,4) -- (0,6);
    \draw (0,4) -- (0,6);
    \draw (2,4) -- (0,6);}

\begin{center}
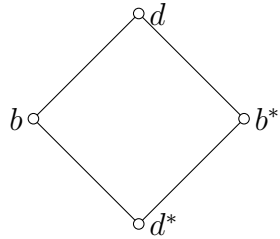

\begin{tikzpicture}[scale=0.7]
    \drawconnectionsa{}
    \draw[black, fill=white] (0,0) circle (.1cm) node[right] {$d^*$};
    \draw[black, fill=white] (-2,2) circle (.1cm) node[left] {$b$};
    \draw[black, fill=white] (2,2) circle (.1cm) node[right] {$b^*$};
    \draw[black, fill=white] (0,4) circle (.1cm) node[right] {$d$};
    \end{tikzpicture}
    \captionof{figure}{Structure of $H:B^2\to D_n$ if $H$ is self-dual.}\label{fig1}
\end{center}

Consider now two elements $b,d\in D_n$, and suppose that $d\ge b|b^*$, see Figure \ref{fig1}.
Then the concatenation 
$$d^*\cdot b\cdot b^*\cdot d$$
is self-dual and represents an element in $D_{n+2}$. Indeed, by Lemma~\ref{L0}, we have
$d^*\le b^*\&b$.  Therefore, we have proved the following lemma.

\begin{lem}\label{L3} The number of self-duals
$\lambda_{n+2}$ is equal to the number of pairs $b,d\in D_n$ that satisfy condition $d\ge b|b^*$. In other words
$$\lambda_{n+2}=\sum_{b\in D_n}re[b|b^*,\top].$$
Furthermore, Lemma~\ref{L2} implies
$$\lambda_{n+2}=\sum_{b\in R_n}\gamma(b)\cdot re[b|b^*,\top].$$
Here we identify each class $[b]\in R_n$ with its representative.
\end{lem}
\section{Counting functions from $B^3$ into $D_n$}
By Lemma~\ref{L1}, the poset $D_{n+3}$ is isomorphic to the poset $D_n^{B^3}$---the set of monotone  functions from 
$B^3=\{000,001,010,100,110,101,011,111\}$ to $D_n$. 

Consider a monotone function $H:B^3\to D_n$. It can be represented as the concatenation:
$$H(000)\cdot H(001)\cdot H(010) \cdot H(011)\cdot H(100)\cdot H(101)\cdot H(110) \cdot H(111)$$
and its dual as
$$H(111)^*\cdot H(110)^*\cdot H(101)^* \cdot H(100 )^*\cdot H(011)^*\cdot H(010)^*\cdot H(001)^* \cdot H(000 )^*.$$
If $H$ is self-dual, then $H(000)=H(111)^*$, $H(001)=H(110)^*$, $H(010)=H(101)^*$ and $H(011)=H(100)^*$.

\begin{lem}\label{L4}
For each $a\in D_n$, such that $a\le a^*$,

for each pair $b,c\in D_n$ such that $a\le b\le c\le a^*$,

for each $d\in D_n$ such that $a\le d\le c\& b^*$,

the concatenation
$$a\cdot b\cdot d\cdot c\cdot c^*\cdot d^*\cdot b^*\cdot a^*$$
represents a self-dual in $D_{n+3}$, see Figure \ref{fig2}.  
\end{lem}

\begin{center}
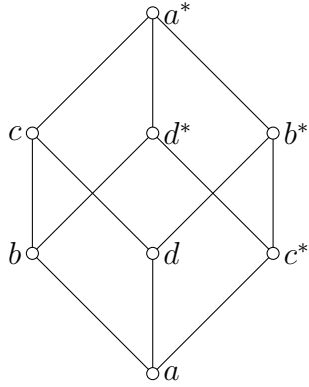

\begin{tikzpicture}[scale=0.8]
\drawconnectionsb{}
\draw[black, fill=white] (0,0) circle (.1cm) node[right] {$a$};
\draw[black, fill=white] (-2,2) circle (.1cm) node[left] {$b$};
\draw[black, fill=white] (0,2) circle (.1cm) node[right] {$d$};
\draw[black, fill=white] (2,2) circle (.1cm) node[right] {$c^*$};
\draw[black, fill=white] (-2,4) circle (.1cm) node[left] {$c$};
\draw[black, fill=white] (0,4) circle (.1cm) node[right] {$d^*$};
\draw[black, fill=white] (2,4) circle (.1cm) node[right] {$b^*$};
\draw[black, fill=white] (0,6) circle (.1cm) node[right] {$a^*$};

\end{tikzpicture}
\captionof{figure}{Structure of $H:B^3\to D_n$ if $H$ is self-dual.}\label{fig2}
\end{center}

\begin{lem}\label{L5} The number of self-duals
$\lambda_{n+3}$ is equal to the number of 4-tuples $(a,b,c,d)\in D_n^4$ satisfying the conditions:
\begin{enumerate}
\item $a\le b\le c\le a^*$
\item $a\le d\le c\& b^*$
\end{enumerate}
\end{lem}

Let
$$G(a)=\sum_{\substack{b,c \in D_n\\a\le b\le c\le a^*}}re[a,c\& b^*].$$
Observe that $G(a)$ is the number of self-dual functions $H\in D_n^{B^3}$, with $H(000)=a$ and $H(111)=a^*$.
$$\lambda_{n+3}=\sum_{\substack{a\in D_n\\a\le a^*}}G(a).$$
Furthermore, Lemma~\ref{L2} implies that for any two equivalent elements $a_1\sim  a_2$ we have $G(a_1)=G(a_2)$.
Hence,   we have
$$\lambda_{n+3}=\sum_{\substack{a\in R_n\\a\le a^*}}\gamma(a)\cdot G(a).$$
Here again we identify the class $[a]\in R_n$ with its representative.

Observe that $a\le a^*$, only if $\ell(a)\le 2^{n-1}$.
Furthermore, if $\ell(a)=2^{n-1}$ and $a\le a^*$, then $a=a^*$, and all
elements $b,c,d,b^*,c^*,d^*$ are equal to $a$; and we have only one self-dual function
$H\in D_n^{B^3}$ with $H(000)=H(111)=a$. Hence,

$$\lambda_{n+3}=\lambda_n+\sum_{\substack{a\in R_n\\\ell(a)<2^{n-1}\\a\le a^*}}\gamma(a)\cdot G(a).$$

\begin{lem}\label{L6}
For each $a\in D_n$, such that $a\le a^*$,

for each $d\in D_n$ such that $a\le d\le a^*$,

for each $b,c\in D_n$ such that: $a\le b\le c\le a^*$,
$b\le d^*$, and $c\ge d$,
the concatenation
$$a\cdot b\cdot \cdot d\cdot c\cdot c^*\cdot d^*\cdot b^*\cdot a^*$$
represents a self-dual in $D_{n+3}$, see Figure \ref{fig2}.  
\end{lem}

\section{Counting functions from $B^4$ into $D_n$}
By Lemma~\ref{L1}, the poset $D_{n+4}$ is isomorphic to the poset $D_n^{B^4}$---the set of monotone  functions from

$B^4=\{0000,0001,0010,0011, 0100, 0101,0110,0111,  1000,1001,1010,1011,$

$   1100,1101,1110,1111\}$  to $D_n$.

Consider a monotone function $H:B^4\to D_n$. It can be represented as the concatenation:
$$H(0000)\cdot H(0001)\cdot H(0010) \cdot H(0011)\cdot H(0100)\cdot H(0101)\cdot H(0110) \cdot H(0111)\cdot$$
$$H(1000)\cdot H(1001)\cdot H(1010) \cdot H(1011)\cdot H(1100)\cdot H(1101)\cdot H(1110) \cdot H(1111)$$

and its dual as
$$H(1111)^*\cdot H(1110)^*\cdot H(1101)^* \cdot H(1100 )^*\cdot H(1011)^*\cdot H(1010)^*\cdot H(1001)^* \cdot H(1000 )^*\cdot$$
$$H(0111)^*\cdot H(0110)^*\cdot H(0101)^* \cdot H(0100 )^*\cdot H(0011)^*\cdot H(0010)^*\cdot H(0001)^* \cdot H(0000 )^*.$$

If $H$ is self-dual then
$H(0000)=H(1111)^*$,
$H(0001)=H(1110)^*$,
$H(0010)=H(1101)^*$,
$H(0011)=H(1100)^*$,
$H(0100)=H(1011)^*$, 
$H(0101)=H(1010)^*$, 
$H(0110)=H(1001)^*$,
$H(0111)=H(1000)^*$.

\begin{lem}\label{L41}
For each $a,b,c\in D_n$, 

for each $h\in D_n$ such that $h\ge a|b|c|a^*|b^*|c^*$,

for each $d,e,f,g\in D_n$ such that 

$a|b|c\le d\le h$

$a|b^*|c^*\le e\le h$

$b|a^*|c^*\le f\le h$

$c|a^*|b^*\le g\le h$

the concatenation
$$h^*\cdot g^*\cdot f^*\cdot a\cdot e^*\cdot b\cdot c\cdot d\cdot d^*\cdot c^*\cdot b^*\cdot e\cdot a^*\cdot f\cdot g\cdot h$$
represents a self-dual in $D_{n+4}$, see Figure \ref{fig3}.  
\end{lem}

\def\drawconnections#1#2#3{\draw (#1,#2) -- (#1-2,#2+2);
    \draw (#1,#2) -- (#1,#2+2);
    \draw (#1,#2) -- (#1+2,#2+2);
    \draw (#1-2,#2+2) -- (#1-2,#2+4);
    \draw (#1-2,#2+2) -- (#1,#2+4);
    \draw (#1,#2+2) -- (#1-2,#2+4);
    \draw (#1,#2+2) -- (#1+2,#2+4);
    \draw (#1+2,#2+2) -- (#1,#2+4);
    \draw (#1+2,#2+2) -- (#1+2,#2+4);
    \draw (#1-2,#2+4) -- (#1,#2+6);
    \draw (#1,#2+4) -- (#1,#2+6);
    \draw (#1+2,#2+4) -- (#1,#2+6);}

\def\drawconnectionsd{
\draw (0,0) -- (5,2);
\draw (-2,2) -- (3,4);
\draw (0,2) -- (5,4);
\draw (2,2) -- (7,4);
\draw (-2,4) -- (3,6);
\draw (0,4) -- (5,6);
\draw (2,4) -- (7,6);
\draw (0,6) -- (5,8);}

\begin{center}
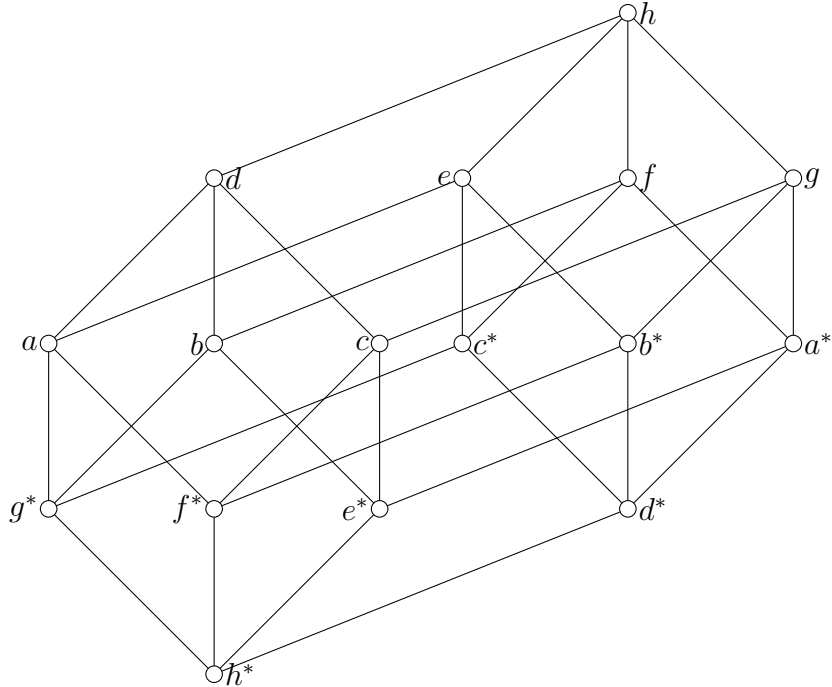

\begin{tikzpicture}[scale=1.1]
\drawconnections{0}{0}{0}
\drawconnectionsd
\draw[black, fill=white] (0,0) circle (.1cm) node[right] {$h^*$};
\draw[black, fill=white] (-2,2) circle (.1cm) node[left] {$g^*$};
\draw[black, fill=white] (0,2) circle (.1cm) node[left] {$f^*$};
\draw[black, fill=white] (2,2) circle (.1cm) node[left] {$e^*$};
\draw[black, fill=white] (-2,4) circle (.1cm) node[left] {$a$};
\draw[black, fill=white] (0,4) circle (.1cm) node[left] {$b$};
\draw[black, fill=white] (2,4) circle (.1cm) node[left] {$c$};
\draw[black, fill=white] (0,6) circle (.1cm) node[right] {$d$};

\drawconnections{5}{2}{0}
\draw[black, fill=white] (5,2) circle (.1cm) node[right] {$d^*$};
\draw[black, fill=white] (3,4) circle (.1cm) node[right] {$c^*$};
\draw[black, fill=white] (5,4) circle (.1cm) node[right] {$b^*$};
\draw[black, fill=white] (7,4) circle (.1cm) node[right] {$a^*$};
\draw[black, fill=white] (3,6) circle (.1cm) node[left] {$e$};
\draw[black, fill=white] (5,6) circle (.1cm) node[right] {$f$};
\draw[black, fill=white] (7,6) circle (.1cm) node[right] {$g$};
\draw[black, fill=white] (5,8) circle (.1cm) node[right] {$h$};

\end{tikzpicture}
\captionof{figure}{Structure of $H:B^4\to D_n$ if $H$ is self-dual.}\label{fig3}
\end{center}

\begin{lem}\label{L51} The number of self-duals
$$\lambda_{n+4}=\sum_{a,b,c\in D_n}\sum_{\substack{h\in D_n \\ h \ge (a|b|c|a^*|b^*|c^*)}} re[a|b|c,h]\cdot re[a|b^*|c^*,h]\cdot re[b|a^*|c^*,h]\cdot re[c|a^*|b^*,h].$$

\end{lem}

\begin{lem}\label{L42}

For each $h\in D_n$ such that $h\ge h^*$,

for each $a,b,c\in D_n$,  $h^*\le a,b,c \le h$

for each $d,e,f,g\in D_n$ such that 

$a|b|c\le d\le h$

$a|b^*|c^*\le e\le h$

$b|a^*|c^*\le f\le h$

$c|a^*|b^*\le g\le h$

the concatenation
$$h^*\cdot g^*\cdot f^*\cdot a\cdot e^*\cdot b\cdot c\cdot d\cdot d^*\cdot c^*\cdot b^*\cdot e\cdot a^*\cdot f\cdot g\cdot h$$
represents a self-dual in $D_{n+4}$; see Figure \ref{fig3}.  
\end{lem}

Let
$$F(h)=\sum_{\substack{h\in D_n \\ h \ge h^*}} \sum_{\substack{a, b, c \in D_n \\ h \ge a,b,c \ge h^*}} re[a|b|c,h]\cdot re[a|b^*|c^*,h]\cdot re[b|a^*|c^*,h]\cdot re[c|a^*|b^*,h].$$
Observe that $F(h)$ is the number of self-dual functions $H\in D_n^{B^4}$, with $H(0000)=h^*$ and $H(1111)=h$.
$$\lambda_{n+4}=\sum_{\substack{h\in D_n\\h^*\le h}}F(h)$$
Furthermore, Lemma~\ref{L2} implies that for any two elements $h_1\sim  h_2$ we have $F(h_1)=F(h_2)$.
Hence,  we have
$$\lambda_{n+4}=\sum_{\substack {h\in R_n\\h\ge h^*}}\gamma(h)\cdot F(h).$$
Here again we identify the class $[h]\in R_n$ with its representative.
Observe that $h\ge h^*$, only if $\ell(h)\ge 2^{n-1}$.
Furthermore, if $\ell(h)=2^{n-1}$ and $h\ge h^*$, then $h=h^*$, and
we have only one self-dual function
$H\in D_n^{B^4}$ with $H(0000)=H(1111)=h$. Hence,

\begin{equation}\label{alg3}
  \lambda_{n+4}=\lambda_n+\sum_{\substack {h\in R_n\\h^*\le h\\\ell(h)> 2^{n-1}}}\gamma\cdot F(h)  
\end{equation}

\section{Testing of algorithms}

In this section we present three algorithms based on the results from the previous sections. We implemented the algorithms in Rust and ran them on a 32-thread Xeon CPU.

\begin{algorithm}[H]
\caption{Calculation of $\lambda_{n+2}$}
    \vspace{1mm}
    \hspace*{\algorithmicindent} \textbf{Input:} $R_n$ with $re[x, \top]$ for all $x \in R_n$ \\
    \hspace*{\algorithmicindent} \textbf{Output:} $s = \lambda_{n+2}$
\begin{algorithmic}[1]
\State Initialize $s = 0$,
\ForAll {$b \in R_n $}
\State $s = s + re[b|b^*,\top] \cdot \gamma(b)$
\EndFor
\end{algorithmic}
\end{algorithm}

Algorithm 1 is based on Lemma \ref{L3}. After loading the preprocessed data into main memory, $\lambda_9$ was computed in 15 seconds. However, preprocessing (the calculation of $R_7$ and its intervals) took approximately 2,5 hours.

\begin{algorithm}[H]
\caption{Calculation of $\lambda_{n+4}$}
    \vspace{1mm}
    \hspace*{\algorithmicindent} {\textbf{Input:} $D_n$; $R_n$; $re[x, y]$ for all $(x, y) \in D_n \times D_n$} \\
    \hspace*{\algorithmicindent} \textbf{Output:} $s = \lambda_{n+4}$
\begin{algorithmic}[1]
\State Initialize $s = 0$,
\ForAll {$a \in R_n $}
\ForAll {$b \in D_n $}
\ForAll {$c \in D_n $}
\ForAll {$h \in D_n, h \ge (a | b | c | a^* | b^* | c^*)$}
\State $s = s + re[a|b|c,h]\cdot re[a|b^*|c^*,h]\cdot re[b|a^*|c^*,h]\cdot re[c|a^*|b^*,h] \cdot \gamma(a)$
\EndFor
\EndFor
\EndFor
\EndFor
\end{algorithmic}
\end{algorithm}

Algorithm 2 is based on Lemma \ref{L51}. Using our implementation of the algorithm, we calculated $\lambda_9$ in 76 seconds, and the preprocessing was almost instantaneous.

\begin{algorithm}[H]
\caption{Calculation of $\lambda_{n+4}$}
    \vspace{1mm}
    \hspace*{\algorithmicindent} {\textbf{Input:} $D_n$; $R_n$; $re[x, y]$ for all $(x, y) \in D_n \times D_n$} \\
    \hspace*{\algorithmicindent} \textbf{Output:} $s = \lambda_{n+4}$
\begin{algorithmic}[1]
\State Initialize $s = \lambda_n$,
\ForAll {$h \in R_n, h^*\le h, \ell(h)> 2^{n-1} $}
\ForAll {$a \in D_n, a < h$}
\ForAll {$b \in D_n, b < h$}
\ForAll {$c \in D_n, c < h$}
\State $s = s + re[a|b|c,h]\cdot re[a|b^*|c^*,h]\cdot re[b|a^*|c^*,h]\cdot re[c|a^*|b^*,h] \cdot \gamma(h)$
\EndFor
\EndFor
\EndFor
\EndFor
\end{algorithmic}
\end{algorithm}

Algorithm 3 is based on Equation \ref{alg3}. The calculation of $\lambda_9$ using our implementation of the algorithm lasted approximately 25 minutes. 

In all cases, we have obtained the following value:

$$ \lambda_9 = 423295099074735261880,$$

which confirms the result of Brouwer et al. \cite{brouwer}.


    
    
    
    
    
    

\end{document}